\documentclass[final,tims]{elsarticle}
\usepackage{amsfonts}
\usepackage{graphicx}
\usepackage{amsmath}
\usepackage{amssymb}
\usepackage{latexsym}
\usepackage{mathrsfs}
\usepackage[all]{xy}
\usepackage{fullpage}
\usepackage{mathtools}
\newtheorem{Theorem}{Theorem}
\newtheorem{cor}{Corollary}

\newtheorem{prop}{Proposition}

\begin{document}
\begin{frontmatter}

\title{Exponential-type Inequalities Involving Ratios of the Modified Bessel Function of the First Kind and their Applications} 

\author[label1]{Prakash Balachandran\corref{cor1}}
\ead{prakashb@math.bu.edu}
\author[label1]{Weston Viles}
\ead{wesviles@bu.edu}
\author[label1]{Eric D. Kolaczyk}
\ead{kolaczyk@bu.edu}
\cortext[cor1]{Corresponding author}

\address[label1]{Boston University, Department of Mathematics and Statistics, 111 Cummington Mall, Boston, MA 02215}

\begin{abstract}
The modified Bessel function of the first kind, $I_{\nu}(x)$, arises in numerous areas of study, such as physics, signal processing, probability, statistics, etc.  As such, there has been much interest in recent years in deducing properties of functionals involving $I_{\nu}(x)$, in particular, of the ratio ${I_{\nu+1}(x)}/{I_{\nu}(x)}$, when $\nu,x\geq 0$.  In this paper we establish sharp upper and lower bounds on $H(\nu,x)=\sum_{k=1}^{\infty} {I_{\nu+k}(x)}/{I_\nu(x)}$ for $\nu,x\geq 0$ that appears as the complementary cumulative hazard function for a Skellam$(\lambda,\lambda)$ probability distribution in the statistical analysis of networks.  Our technique relies on bounding existing estimates of ${I_{\nu+1}(x)}/{I_{\nu}(x)}$ from above and below by quantities with nicer algebraic properties, namely exponentials, to better evaluate the sum, while optimizing their rates in the regime when $\nu+1\leq x$ in order to maintain their precision.  We demonstrate the relevance of our results through applications, providing an improvement for the well-known asymptotic $\exp(-x)I_{\nu}(x)\sim {1}/{\sqrt{2\pi x}}$ as $x\rightarrow \infty$, upper and lower bounding $\mathbb{P}\left[W=\nu\right]$ for $W\sim Skellam(\lambda_1,\lambda_2)$,  and deriving a novel concentration inequality on the $Skellam(\lambda,\lambda)$ probability distribution from above and below.
\end{abstract}


\end{frontmatter}

\section{Introduction} 
\label{Introduction}

The modified Bessel function of the first kind, $$I_{\nu}(x)=\sum_{k=0}^{\infty} \frac{1}{k! \Gamma(k+\nu+1)} \left(\frac{x}{2}\right)^{2k+\nu},$$  arises in numerous applications.  In elasticity \cite{Simpson}, one is interested in ${I_{\nu+1}(x)}/{I_{\nu}(x)}$.  In image-noise modeling \cite{Hwang}, denoising photon-limited image data \cite{Wolfe}, sports data \cite{Karlis}, and statistical testing \cite{Skellam}, one is interested in $I_{\nu}(x)$ as it arises in the kernel of the probability mass function of a Skellam probability distribution.  The functions $I_0(x)$ and $I_1(x)$ arise as rates in concentration inequalities in the behavior of sums of independent $\mathbb{R}^N$-valued, symmetric random vectors \cite{Kanter} \cite{Cranston}.   Excellent summaries of applications of $I_{\nu}(x)$ in probability and statistics may be found in~\cite{Robert,Johnson}.  For example, these functions are used in the determination of maximum likelihood and minimax estimators for the bounded mean in~ \cite{Marchand1,Marchand2}. Furthermore,  \cite{Yuan} gives applications of the modified Bessel function of the first kind to the so-called Bessel probability distribution, while~\cite{Simon} has applications of $I_{\nu}(x)$ in the generalized Marcum Q-function that arises in communication channels.  Finally, \cite{Fotopoulos} gives applications in finance.

With $I_{\nu}(x)$ regularly arising in new areas of application comes a corresponding need to continue to better understand its properties, both as a function of $\nu$ and of $x$.  There is, of course, already much that has been established on this topic.  For example, \cite{Luke} first studied inequalities on generalized hypergeometric functions and was the first to give upper and lower bounds on $I_{\nu}(x)$ for $x>0$ and $\nu>-{1}/{2}$.  Additionally, \cite{Amos} was concerned with computation of $I_{\nu}(x)$, and provided a way to produce rapid evaluations of ratios $I_{\nu+1}(x)/I_{\nu}(x)$, and hence $I_{\nu}(x)$ itself, through recursion.  Several other  useful representations of $I_{\nu}(x)$ are also provided in~\cite{Amos}.  More recently, various convexity properties of $I_{\nu}(x)$ have been studied in~\cite{Neuman} and \cite{Baricz Neuman}.  In the last decade, motivated by results in finite elasticity, \cite{Simpson} and~\cite{Laforgia} provide bounds on $I_{\nu}(x)/I_{\nu}(y)$ for $\nu>0$ and $0<x<y$, while \cite{Baricz Sum} provides bounds on the quantity $\exp(-x)x^{-\nu} \left[I_{\nu}(x)+I_{\nu+1}(x)\right]$ arising in concentration of random vectors, as in~\cite{Kanter, Cranston}.  Motivated by applications in communication channels, \cite{Baricz Marcum Q} develops bounds on the generalized Marcum-Q function, and the same author in \cite{Baricz Turanians} develops estimates on the so-called Turan-type inequalities $I^2_{\nu}(x)-I_{\nu-1}(x)I_{\nu+1}(x)$.  For an excellent review of modern results on $I_{\nu}(x)$ and its counterpart $K_{\nu}(x)$, the modified Bessel function of the second kind, we refer the reader to \cite{Baricz Edin}.

In our own ongoing work in the statistical analysis of networks, the function $I_{\nu}(x)$ has arisen as well, in a manner that -- to the best of our knowledge -- has yet to be encountered and addressed in the literature.  Specifically, in seeking to establish the probability distribution of the discrepancy between (a) the true number of edges in a network graph, and (b) the number of edges in a `noisy' version of that graph, one is faced with the task of analyzing the distribution of the difference of two sums of dependent binary random variables.  Under a certain asymptotic regime, it is reasonable to expect that each sum converge to a certain Poisson random variable and, hence, their difference, to a so-called Skellam distribution.  The latter is the name for the probability distribution characterized by the difference of two independent Poisson random variables and --  notably -- has a kernel defined in terms of $I_{\nu}(x)$ \cite{Skellam} for $x\geq 0$ and $\nu \in \mathbb{N}$.  One way to study the limiting behavior of our difference of sums is through Stein's method \cite{Barbour}.  As part of such an analysis, however, non-asymptotic upper bounds are necessary on the quantity 
\begin{equation}
\label{H Def}
H(\nu, x)=\sum_{k=1}^{\infty} \frac{I_{\nu+k}(x)}{I_{\nu}(x)} \enskip ,
\end{equation}
for $\nu \in \mathbb{N}$ that have a scaling of $\sqrt{x}$ for $\nu$ near 0.  Unfortunately, using current bounds on $I_{\nu+1}(x)/I_{\nu}(x)$ to lower and upper bound the infinite sum in $H(\nu,x)$ in (\ref{H Def}) for $\nu,x\geq 0$ necessitate the use of a geometric series-type argument, the resulting expressions of which both do not have  this kind of behavior near $\nu=0$.  In particular, we show that such an approach, for $\nu=0$, yields a lower bound that is order one and an upper bound that is order $x$ as $x\rightarrow \infty$.  See (\ref{Geometric Bounds}) below.

The purpose of this paper is to derive bounds on $I_{\nu+1}(x)/I_{\nu}(x)$ which, when used to lower and upper bound the infinite sum arising in $H(\nu,x)$, lead to better estimates on $H(\nu,x)$ near $\nu=0$ compared to those obtained using current estimates, for $\nu,x\geq 0$ and $\mathbb{\nu}\in \mathbb{R}$.  In particular, we show that it is possible to derive both upper and lower bounds on $H(\nu,x)$ that behave as $\sqrt{x}$ for $x$ large.  When we restrict $\nu$ to $\mathbb{N}$, we can apply these results to obtain a concentration inequality for the Skellam distribution, to bound the probability mass function of the Skellam distribution, and to upper and lower bound $\exp(-x)I_{\nu}(x)$ for any $\nu,x\geq 0$, improving on the asymptotic $\exp(-x)I_{\nu}(x)\sim {1}/{\sqrt{2\pi x}}$ as $x\rightarrow \infty$ in \cite{AS}, at least for $\nu \in \mathbb{N}$.

In our approach to analyzing the function $H(\nu,x)=\sum_{n=1}^{\infty} {I_{\nu+n}(x)}/{I_{\nu}(x)}$, we first write each term in the sum using the iterative product,
\begin{equation}
\label{intro: iterative product}
\frac{I_{\nu+n}(x)}{I_{\nu}(x)}=\prod_{k=0}^{n-1} \frac{I_{\nu+k+1}(x)}{I_{\nu+k}(x)}
\end{equation}
and split the infinite sum (\ref{H Def}) into two regimes: one where $[\nu]+2>x$ and the other when $[\nu]+2\leq x$, where $[\nu]$ denotes the floor function of $\nu$.  In the former regime, the "tail'' of $H(\nu,x)$, we can use existing estimates on $I_{\nu+1}(x)/I_{\nu}(x)$ in a geometric series to lower and upper bound $H(\nu,x)$ in a way that preserves the scaling of $H(\nu,x)$ in $\nu$ and $x$.  In the latter regime, lower and upper bounds on the function ${I_{\nu+1}(x)}/{I_{\nu}(x)}$ for $\nu \in \mathbb{R}$ and $\nu,x\geq 0$ are now required with algebraic properties suitable to better sum the the products (\ref{intro: iterative product}) arising in $H(\nu,x)$ in a way that preserves the behavior of $H(\nu,x)$ near $\nu=0$ for large $x$.  

To provide these bounds on $I_{\nu+1}(x)/I_{\nu}(x)$, we begin with those in \cite{Amos}, valid for  $\nu,x\geq 0$, which can be expressed as
\begin{equation}
\label{Amos Bound}
\sqrt{1+\left(\frac{\nu+1}{x}\right)^2} -\frac{\nu+1}{x} \leq  \frac{I_{\nu+1}(x)}{I_{\nu}(x)}\leq  \sqrt{1+\left(\frac{\nu+\frac{1}{2}}{x}\right)^2} -\frac{\nu+\frac{1}{2}}{x} \enskip ,
\end{equation}
and weaken them to those with nicer, exponential properties, using a general result on the best exponential approximation for the function $f(x)=\sqrt{1+x^2}-x$  for $x\in [0,1]$.  When applied to $I_{\nu+1}(x)/I_{\nu}(x)$ for $\nu+1\leq x$, we obtain

$$\exp\left(-\frac{\nu+1}{x}\right)\leq \frac{I_{\nu+1}(x)}{I_{\nu}(x)} \leq \exp\left(-\alpha_0 \frac{\nu+\frac{1}{2}}{x}\right).$$
See Proposition~\ref{besselexp} and Corollary~\ref{Bessel Ratio}.

Using these bounds to lower and upper bound $H(\nu,x)$ described in the above fashion, we obtain

\begin{enumerate}
\item  For any $\nu,x\geq 0$ (and in particular, for $[\nu]+2>[x]$),

\begin{equation}
\label{Geometric Bounds}
F(\nu+1,x)(1+F(\nu+2,x))\leq H(\nu,x) \leq \frac{F\left(\nu+\frac{1}{2},x\right)}{1-F\left(\nu+\frac{3}{2},x\right) }.
\end{equation}

\item  If $\nu,x\geq 0$ and $[\nu]+2\leq [x]$, 

\begin{equation}
\label{Our H Bounds}
\mathcal{L}(\nu,x)\leq H(\nu,x)\leq \mathcal{U}(\nu,x)
\end{equation}

where 

\begin{equation}
\label{Our Lower H Bound}
\begin{aligned}
\mathcal{L}(\nu,x)&=   \frac{2xe^{-\frac{1}{x}\left(\nu+1\right)}}{\nu+\frac{3}{2}+\sqrt{\left(\nu+\frac{3}{2}\right)^2+4x}}- \frac{2xe^{-\frac{1}{2x} ([x]-\nu-\nu_f+1)([x]+\nu-\nu_f+2)}}{[x]+\frac{3}{2}+\sqrt{ \left([x]+\frac{3}{2}\right)^2+\frac{8x}{\pi}}}\\
&\hspace{2in}+e^{-\frac{([x]-[\nu]-1)([x]+\nu+\nu_f)}{2x}}F([x]+\nu_f,x)(1+F([x]+\nu_f+1,x))\enskip.
\end{aligned}
\end{equation}

and

\begin{equation}
\label{Our Upper H Bound}
\begin{aligned}
\mathcal{U}(\nu,x)&= \frac{2x}{\alpha_0}\left[\frac{1}{\nu +\sqrt{\nu^2+\frac{8x}{\pi\alpha_0}}}-\frac{e^{-\frac{\alpha_0}{2x} (x^2-\nu^2)}}{x+\sqrt{x^2+\frac{4x}{\alpha_0}}}\right] \\
&\hspace{1in}+e^{-\frac{\alpha_0}{2x} ([x]-[\nu]-1)([x]+\nu+\nu_f-1)}\frac{F\left([x]+\nu_f-\frac{1}{2},x\right)}{1-F\left([x]+\nu_f+\frac{1}{2},x\right)}\enskip.
\end{aligned}
\end{equation}

\end{enumerate}

where $\alpha_0=-\log(\sqrt{2}-1)$, $[x]$ denotes the floor function of $x$, $x=[x]+x_f$, and $$F(\nu,x)=\frac{x}{\nu+\sqrt{\nu^2+x^2}}\enskip.$$

We note here that the bounds in (\ref{Our Lower H Bound}) and (\ref{Our Upper H Bound}) are similar to those occurring in (\ref{Geometric Bounds}), but now with exponentially decaying factors in $x$ plus an incurred error from $\sum_{n=1}^{[x]-[\nu]-1} \prod_{k=0}^{n-1} {I_{\nu+k+1}(x)}/{I_{\nu+k}(x)}$ which behaves like a partial sum of a Gaussian over the integers from $\nu$ to $[x]$.  These contributions are lower and upper bounded by the first differences in both (\ref{Our Lower H Bound}) and (\ref{Our Upper H Bound}), and are responsible for our bounds behaving like $\sqrt{x}$ for $\nu$ near $0$ and $x$ large.  Indeed, if one were to simply use (\ref{Geometric Bounds}) for all $\nu,x\geq 0$, then for $\nu=0$ and large $x$, the lower bound is order 1, while the upper bound is of order $x$.

The rest of this paper is organized as follows.  We derive our bounds in Section~\ref{main results},  give applications in Section~\ref{appls}, and provide some discussion in Section~\ref{conclusion}.   In Section~\ref{Bessel function}, we first give the result on the best exponential approximation to the function $f(x)=\sqrt{1+x^2}-x$ for $x\in [0,1]$,  in Proposition~\ref{besselexp}, and apply them to lower and upper bounding $I_{\nu+1}(x)/I_{\nu}(x)$ when $\nu+1\leq x$, obtaining Corollary~\ref{Bessel Ratio}.  In Section~\ref{Hazard function}, we then use these bounds to give the upper and lower bounds on $H(\nu,x)$ for $\nu,x\geq 0$.  Combining these bounds with a normalizing condition from the Skellam distribution, we provide in Section~\ref{Asymptotic application} deterministic upper and lower bounds on $\exp(-x)I_{\nu}(x)$ for $\nu\in \mathbb{N}$ and apply them to obtain upper and lower bounds on $\mathbb{P}\left[W=\nu\right]$ for $W\sim Skellam(\lambda_1,\lambda_2)$. Finally, in Section~\ref{Skellam application}, we apply the results on $H(\nu,x)$ to deriving a concentration inequality for the $Skellam(\lambda,\lambda)$.

\section{Main Results: Bounds}
\label{main results}

\subsection{Pointwise bounds on ${I_{\nu+1}(x)}/{I_{\nu}(x)}$}
\label{Bessel function}
We  begin with upper and lower bounds on the ratio ${I_{\nu+1}(x)}/{I_{\nu}(x)}$.  First, we need the following Proposition.

\begin{prop}{(Best Exponential Approximation)}
\label{besselexp}

For all $x\in [0,1]$,
\begin{equation}
\label{BEA 1}
\exp(-x)\leq  \sqrt{1+x^2}-x \leq \exp(-\alpha_0 x) \enskip,
\end{equation}

where $\alpha_0=-\log(\sqrt{2}-1)\approx 0.8814.$  Moreover, these are the best possible arguments of the exponential, keeping constants of $1$.  
\end{prop}

{\bf Proof of Proposition \ref{besselexp}:}

We want to find the best constants $\alpha_1,\alpha_2>0$ for which $$\exp(-\alpha_1 x)\leq \sqrt{1+x^2}-x \leq \exp(-\alpha_2 x), \; \; \; x\in [0,1].$$

To this end, consider the function $$f(x)=\left(\sqrt{1+x^2}-x\right) \exp(\alpha x)$$ for some $\alpha>0$.  We want to find the maximum and minimum values of $f(x)$ on the interval $[0,1]$.  First, note that $$f(0)=1$$ and $$f(1)=(\sqrt{2}-1)\exp(\alpha).$$

To check for critical points, we have

$$
\begin{aligned}
f'(x)&=\exp(\alpha x) \left[ \alpha \left(\sqrt{1+x^2}-x\right) +\frac{x}{\sqrt{1+x^2}} -1\right]\\
&=\frac{\exp(\alpha x)}{\sqrt{1+x^2}} \left[ (\alpha+x+\alpha x^2) - (1+\alpha x)\sqrt{1+x^2}\right]\\
&=\frac{\exp(\alpha x)}{\sqrt{1+x^2} \left[(\alpha+x+\alpha x^2)+(1+\alpha x)\sqrt{1+x^2}\right]}\left[ (\alpha+x+\alpha x^2)^2 - (1+\alpha x)^2 (1+x^2)\right]\\
&=0
\end{aligned}
$$

Thus, we require

$$(\alpha+x+\alpha x^2)^2 = (1+\alpha x)^2 (1+x^2).$$

Expanding both sides of this equation and after some algebra, we get

$$\alpha^2 +\alpha^2x^2=1 \Leftrightarrow x=x_0:=\pm \sqrt{\frac{1-\alpha^2}{\alpha^2}}.$$

Furthermore, this computation shows that this value is always a local minimum.

{\bf Case 1: $\alpha\geq 1$}

In this case, there are no critical points, and the function $f(x)$ is monotone increasing on $(0,1)$.  We find that the upper bound is $f(1)=(\sqrt{2}-1)\exp(\alpha)$ and the lower bound is $f(0)=1$, so that

$$\exp(-\alpha x) \leq \sqrt{1+x^2}-x \leq (\sqrt{2}-1)\exp(\alpha) \exp(-\alpha x).$$

The lower bound maximizes at the value $\alpha=1$.

\subsection*{{\bf Case 2:  $\alpha\leq {1}/{\sqrt{2}}$:}}
In this regime, $\alpha={1}/{\sqrt{2}}$, $x_0\geq 1$, and now the function $f$ is monotone decreasing.  Thus,

$$(\sqrt{2}-1)\exp(\alpha) \exp(-\alpha x)\leq \sqrt{1+x^2}-x \leq \exp(-\alpha x)$$

We can minimize the upper bound by taking $\alpha={1}/{\sqrt{2}}$.

\subsection*{{\bf Case 3: ${1}/{\sqrt{2}}\leq \alpha \leq 1$}}

$$f_{x_0}(\alpha)=(1-\sqrt{1-\alpha^2})\frac{\exp\left(\sqrt{1-\alpha^2}\right)}{\alpha}=\frac{\alpha}{1+\sqrt{1-\alpha^2}} \exp(\sqrt{1-\alpha^2}).$$
$$
\begin{aligned}
f_{x_0}'(\alpha)=\frac{\exp\left(\sqrt{1-\alpha^2}\right)(1-\alpha^2)}{(\sqrt{1-\alpha^2})(1+\sqrt{1-\alpha^2})}
\end{aligned}
$$

Thus, we find that starting at $\alpha={1}/{\sqrt{2}}$ the critical point occurs at $x=1$, and monotonically moves to the left at which point it settles at $x=0$ at $\alpha=1$.  While it does this, the value of the local minimum, $f(x_0)$, increases monotonically, as does $f(1)$.

So, in all cases, $ f(0)\geq f(x_0)$ and $f(x_0)\leq  f(1)$.  But, $f(0)\geq f(1)$ for ${1}/{\sqrt{2}}\leq \alpha \leq \alpha_0$ and then $f(0)\leq  f(1)$ for $\alpha_0\leq \alpha \leq 1$ and equality only occurs at $\alpha=\alpha_0$.  Since we are interested in constants of $1$, in the former case, ${1}/{\sqrt{2}}\leq \alpha \leq \alpha_0$ implies,

$$\sqrt{1+x^2}-x \leq \exp(-\alpha x).$$

We can minimize the upper bound by taking $\alpha=\alpha_0$.

Thus,

$$\exp(-x)\leq  \sqrt{1+x^2}-x \leq \exp(-\alpha_0 x)$$ 

where $\alpha_0=-\log(\sqrt{2}-1)\approx 0.8814.$

\begin{flushright}
$\square$
\end{flushright}

Next, applying Proposition \ref{besselexp} to the ratio ${I_{\nu+1}(x)}/{I_{\nu}(x)}$, we have the following corollary,

\begin{cor}
\label{Bessel Ratio}
Let $\nu,x\geq 0$, and let $\alpha_0=-\log(\sqrt{2}-1)$.  If $\nu+1\leq  x$, then

\begin{equation}
\label{Bessel Ratio I}
\exp\left(-\frac{\nu+1}{x}\right)\leq \frac{I_{\nu+1}(x)}{I_{\nu}(x)} \leq \exp\left(-\alpha_0 \frac{\nu+\frac{1}{2}}{x}\right).
\end{equation}
\end{cor}

{\bf Proof of Corollary \ref{Bessel Ratio}:}
Note that by the bounds \cite{Amos}, for $\nu,x\geq 0$,  
\begin{equation}
\label{Amos Bounds}
\begin{aligned}
&\frac{I_{\nu+1}(x)}{I_{\nu}(x)}\geq \sqrt{1+\left(\frac{\nu+1}{x}\right)^2} -\frac{\nu+1}{x} = \frac{x}{\nu+1+\sqrt{x^2+(\nu+1)^2}} \enskip, \\
&\frac{I_{\nu+1}(x)}{I_{\nu}(x)}\leq  \frac{x}{\nu+\frac{1}{2}+\sqrt{x^2+\left(\nu+\frac{1}{2}\right)^2}} =\sqrt{1+\left(\frac{\nu+\frac{1}{2}}{x}\right)^2} -\frac{\nu+\frac{1}{2}}{x}.
\end{aligned}
\end{equation}

We note that we cannot use the more precise lower bound in \cite{Amos}, 

$$
\begin{aligned}
\label{Amos Bounds Optimal}
&\frac{I_{\nu+1}(x)}{I_{\nu}(x)}\geq  \frac{x}{\nu+\frac{1}{2}+\sqrt{x^2+\left(\nu+\frac{3}{2}\right)^2}} 
\end{aligned}
$$
since we require the arguments in $\nu$ in the denominator to be the same.

When $\nu+1\leq  x$, both $(\nu+{1}/{2})/{x}, (\nu+1)/{x}\leq 1$ so that by Proposition \ref{besselexp},

$$\exp\left(-\frac{\nu+1}{x}\right)\leq \frac{I_{\nu+1}(x)}{I_{\nu}(x)} \leq \exp\left(-\alpha_0 \frac{\nu+\frac{1}{2}}{x}\right).$$

\begin{flushright}
$\square$
\end{flushright}

\begin{figure}[htp!]
\includegraphics[width=15cm, height=10cm]{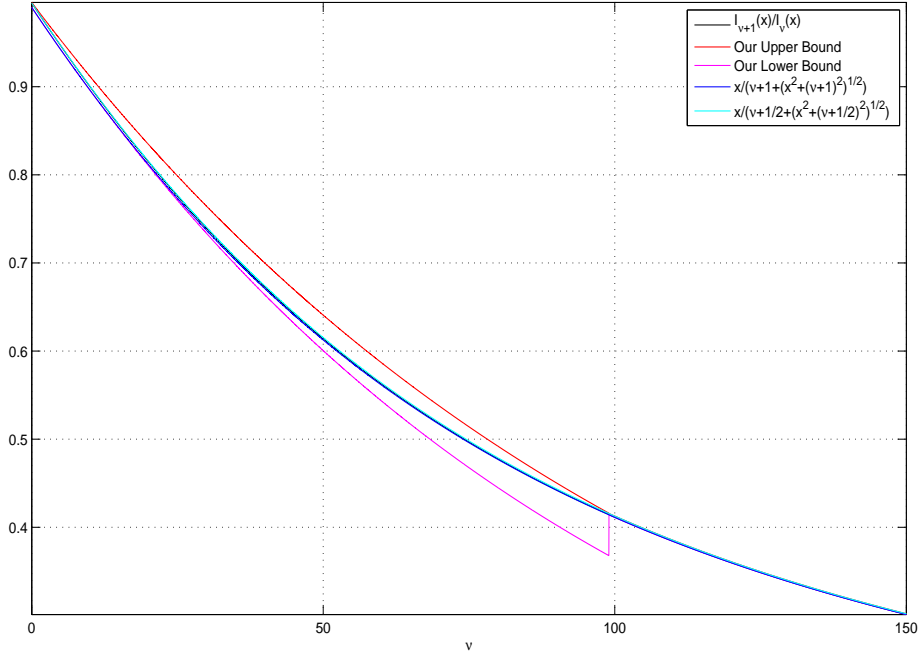}
\caption{An illustration of the exponential-type bounds from Corollary \ref{Bessel Ratio} on ${I_{\nu+1}(x)}/{I_{\nu}(x)}$ for $x=100$, over the interval $[0,150]$ taken in steps of $0.015$.  For $\nu+1\leq x$, we apply the bounds from Corollary \ref{Bessel Ratio}, after which we use the lower and upper bounds ${x}/(\nu+1+\sqrt{x^2+(\nu+1)^2})$ and ${x}/(\nu+\frac{1}{2}+\sqrt{x^2+\left(\nu+\frac{1}{2}\right)^2})$, respectively.  For comparison, we also plot these latter bounds for $\nu+1\leq x$, and note that due to how precise these bounds are for any $\nu$,  the black, blue, and cyan curves nearly coincide.}
\label{besselfunccomp2}
\end{figure}

We illustrate these bounds on ${I_{\nu+1}(x)}/{I_{\nu}(x)}$ in Figure \ref{besselfunccomp2}. 

\newpage

\subsection{Bounds on $H(\nu,x)$}
\label{Hazard function}

The bounds in Section \ref{Bessel function} on $I_{\nu+1}(x)/I_{\nu}(x)$ have extremely nice algebraic properties suitable for evaluation of products.  This allows us to obtain explicit and interpretable bounds on $H(\nu,x)$.  

Recall our program outlined in Section \ref{Introduction}: for any $\nu,x\geq 0$,
$$
\begin{aligned}
H(\nu,x) &= \begin{cases}
\sum_{n=1}^{\infty} \prod_{k=0}^{n-1} \frac{I_{\nu+k+1}(x)}{I_{\nu+k}(x)} & [\nu]+2> [x]\\
\sum_{n=1}^{[x]-[\nu]-1} \prod_{k=0}^{n-1} \frac{I_{\nu+k+1}(x)}{I_{\nu+k}(x)} + \sum_{n=[x]-[\nu]}^{\infty} \prod_{k=0}^{[x]-[\nu]-2} \frac{I_{\nu+k+1}(x)}{I_{\nu+k}(x)} \prod_{k=[x]-[\nu]-1}^{n-1} \frac{I_{\nu+k+1}(x)}{I_{\nu+k}(x)} & [\nu]+2\leq [x] \enskip.
\end{cases}\\
\end{aligned}
$$
Using our bounds on $I_{\nu+1}(x)/I_{\nu}(x)$ in Corollary \ref{Bessel Ratio}, the first term in the regime $[\nu]+2\leq [x]$ behaves like a sum of discrete Gaussians.  The second term in the regime $[\nu]+2\leq [x]$ and the term in the regime $[\nu]+2>x$ are ''tail"-like quantities, and a simple geometric series-type argument using (\ref{Amos Bounds}) suffices to capture the behavior of $H(\nu,x)$.  In fact, such a geometric series-type argument holds for any $\nu\geq 0$, and we give the full result which will be useful for comparison.

\begin{Theorem}
\label{Hboundexp}
Let $\nu,x\geq 0$ and $H(\nu,x)$ be defined as in (\ref{H Def}).  Then, (\ref{Geometric Bounds}), (\ref{Our H Bounds}), (\ref{Our Lower H Bound}) and (\ref{Our Upper H Bound}) hold.  
\end{Theorem}

{\bf Proof of Theorem \ref{Hboundexp}:}

\begin{enumerate}

\item  We first prove (\ref{Geometric Bounds}).  Note that

$$H(\nu,x)=\frac{I_{\nu+1}(x)}{I_{\nu}(x)}(1+H(\nu+1,x)) \enskip,$$

which in view of (\ref{Amos Bounds}) yields,

\begin{equation}
\label{H Recurse}
F(\nu+1,x)(1+H(\nu+1,x))\leq H(\nu,x) \leq F\left(\nu+\frac{1}{2},x\right) (1+H(\nu+1,x))\enskip.
\end{equation}

Next, we have 

$$\begin{aligned}
H(\nu+1,x)&=\sum_{n=1}^{\infty} \prod_{k=1}^{n} \frac{I_{\nu+k+1}(x)}{I_{\nu+k}(x)}\\
&\leq \sum_{n=1}^{\infty} \prod_{k=1}^{n} \frac{x}{\nu+k+\frac{1}{2} + \sqrt{\left(\nu+k+\frac{1}{2}\right)^2+x^2}}\\
&\leq \sum_{n=1}^{\infty} \left( \frac{x}{\nu+\frac{3}{2} + \sqrt{\left(\nu+\frac{3}{2}\right)^2+x^2}}\right)^n\\
&=\sum_{n=1}^{\infty}F\left(\nu+\frac{3}{2},x\right)^n\enskip,
\end{aligned}
$$
so that $$1+H(\nu+1,x)\leq \sum_{n=0}^{\infty}F\left(\nu+\frac{3}{2},x\right)^n=\frac{1}{1-F\left(\nu+\frac{3}{2},x\right)} \enskip.$$

Thus,

$$H(\nu,x) \leq \frac{F\left(\nu+\frac{1}{2},x\right)}{1-F\left(\nu+\frac{3}{2},x\right) }$$

yielding the upper bound in (\ref{Geometric Bounds}).

For the lower bound, note that using (\ref{Amos Bounds}), $$F(\nu+2,x)\leq \frac{I_{\nu+2}(x)}{I_{\nu+1}(x)} \leq H(\nu+1,x)\enskip,$$ which in view of (\ref{H Recurse}) implies

$$F(\nu+1)(1+F(\nu+2,x))\leq H(\nu,x)\enskip.$$

This completes the proof of (\ref{Geometric Bounds}).

\item Next, we prove (\ref{Our H Bounds}), (\ref{Our Lower H Bound}) and (\ref{Our Upper H Bound}).  Note that using an iterated product, we may write

$$H(\nu,x)= \sum_{n=1}^{\infty} \frac{I_{\nu+n}(x)}{I_{\nu}(x)}=\sum_{n=1}^{\infty} \prod_{k=0}^{n-1} \frac{I_{\nu+k+1}(x)}{I_{\nu+k}(x)}$$

so that for $[\nu]+2\leq [x]$,

\begin{equation}
\label{separate H}
\begin{aligned}
H(\nu,x) &=\sum_{n=1}^{[x]-[\nu]-1} \prod_{k=0}^{n-1} \frac{I_{\nu+k+1}(x)}{I_{\nu+k}(x)} + \sum_{n=[x]-[\nu]}^{\infty} \prod_{k=0}^{[x]-[\nu]-2} \frac{I_{\nu+k+1}(x)}{I_{\nu+k}(x)} \prod_{k=[x]-[\nu]-1}^{n-1} \frac{I_{\nu+k+1}(x)}{I_{\nu+k}(x)} \\
&=\sum_{n=1}^{[x]-[\nu]-1} \prod_{k=0}^{n-1} \frac{I_{\nu+k+1}(x)}{I_{\nu+k}(x)} +\prod_{k=0}^{[x]-[\nu]-2} \frac{I_{\nu+k+1}(x)}{I_{\nu+k}(x)} \sum_{n=[x]-[\nu]}^{\infty}  \prod_{k=[x]-[\nu]-1}^{n-1} \frac{I_{\nu+k+1}(x)}{I_{\nu+k}(x)} \enskip.
\end{aligned}
\end{equation}

First, we deal with the sum in the second term.  Using similar arguments as above for the upper bound in the first part of the theorem, we can write

$$
\begin{aligned}
& \sum_{n=[x]-[\nu]}^{\infty} \prod_{k=[x]-[\nu]-1}^{n-1} \frac{I_{\nu+k+1}(x)}{I_{\nu+k}(x)}\\
&\leq\sum_{n=[x]-[\nu]}^{\infty} \prod_{k=[x]-[\nu]-1}^{n-1} \frac{x}{\nu+k+\frac{1}{2}+\sqrt{(\nu+k+\frac{1}{2})^2+x^2}}\\
&=F\left([x]+\nu_f-\frac{1}{2}\right)\cdot \left(1+ \sum_{n=[x]-[\nu]}^{\infty} \prod_{k=[x]-[\nu]}^{n} \frac{x}{\nu+k+\frac{1}{2}+\sqrt{(\nu+k+\frac{1}{2})^2+x^2}}\right)\\ 
&\leq F\left([x]+\nu_f-\frac{1}{2}\right) \cdot \left(1+\sum_{n=[x]-[\nu]}^{\infty} \left(\frac{x}{[x]+\nu_f+\frac{1}{2}+\sqrt{([x]+\nu_f+\frac{1}{2})^2+x^2}}\right)^{n+1-([x]-[\nu])} \right)\\
&=\frac{F\left([x]+\nu_f-\frac{1}{2}\right)}{1-F\left([x]+\nu_f+\frac{1}{2},x\right)}\enskip.
 \end{aligned}
  $$
  
 Similar arguments as in the lower bound for $[\nu]+2>[x]$ yield,
 
 $$
 \begin{aligned}
 &\sum_{n=[x]-[\nu]}^{\infty} \prod_{k=[x]-[\nu]-1}^{n-1} \frac{I_{\nu+k+1}(x)}{I_{\nu+k}(x)}\\
 &\geq \frac{I_{[x]+\nu_f}(x)}{I_{[x]+\nu_f-1}(x)}\left(1+ \frac{I_{[x]+\nu_f+1}(x)}{I_{[x]+\nu_f}(x)}\right)\\
& \geq F([x]+\nu_f,x)(1+F([x]+\nu_f+1,x))\enskip.
 \end{aligned}
 $$
 
 Thus,
 
 $$H(\nu,x)\leq \sum_{n=1}^{[x]-[\nu]-1} \prod_{k=0}^{n-1} \frac{I_{\nu+k+1}(x)}{I_{\nu+k}(x)} +\prod_{k=0}^{[x]-[\nu]-2} \frac{I_{\nu+k+1}(x)}{I_{\nu+k}(x)} \frac{F\left([x]+\nu_f-\frac{1}{2}\right)}{1-F\left([x]+\nu_f+\frac{1}{2},x\right)}$$

and

$$H(\nu,x)\geq \sum_{n=1}^{[x]-[\nu]-1} \prod_{k=0}^{n-1} \frac{I_{\nu+k+1}(x)}{I_{\nu+k}(x)} +\prod_{k=0}^{[x]-[\nu]-2} \frac{I_{\nu+k+1}(x)}{I_{\nu+k}(x)}F([x]+\nu_f,x)(1+F([x]+\nu_f+1,x)) \enskip.$$

Next, note that, each term in each of the products above have $\nu+k+1\leq x$, since the largest $k$ can be is $k=[x]-[\nu]-2$ and $\nu+([x]-[\nu]-2)+1=[x]+\nu_f-1\leq x$.  Thus, we may apply Corollary \ref{Hboundexp} to obtain 

$$
\begin{aligned}
\prod_{k=0}^{[x]-[\nu]-2} \frac{I_{\nu+k+1}(x)}{I_{\nu+k}(x)} &\leq \prod_{k=0}^{[x]-[\nu]-2} e^{-\alpha_0 \frac{\nu+k+\frac{1}{2}}{x}}\\
&=\exp\left( -\frac{\alpha_0}{x} \left( \left(\nu+\frac{1}{2}\right) \left([x]-[\nu]-1\right) + \frac{\left([x]-[\nu]-1\right)\left([x]-[\nu]-2\right)}{2}\right)\right)\\
&=\exp\left( -\frac{\alpha_0}{2x} \left([x]-[\nu]-1\right)\left(2\nu+1+[x]-[\nu]-2\right) \right)\\
&=\exp\left(-\frac{\alpha_0}{2x} ([x]-[\nu]-1)([x]+\nu+\nu_f-1)\right) \enskip.
\end{aligned}
$$

Likewise,
$$
\begin{aligned}
\prod_{k=0}^{[x]-[\nu]-2} \frac{I_{\nu+k+1}(x)}{I_{\nu+k}(x)} &\geq \prod_{k=0}^{[x]-[\nu]-2} e^{- \frac{\nu+k+1}{x}}=\exp\left(-\frac{([x]-[\nu]-1)([x]+\nu+\nu_f)}{2x}\right)\\
\end{aligned}
$$

implying

 $$H(\nu,x)\leq \sum_{n=1}^{[x]-[\nu]-1} \prod_{k=0}^{n-1} \frac{I_{\nu+k+1}(x)}{I_{\nu+k}(x)} +  e^{-\frac{\alpha_0}{2x} ([x]-[\nu]-1)([x]+\nu+\nu_f-1)}\frac{F\left([x]+\nu_f-\frac{1}{2}\right)}{1-F\left([x]+\nu_f+\frac{1}{2},x\right)}$$

and

$$H(\nu,x)\geq \sum_{n=1}^{[x]-[\nu]-1} \prod_{k=0}^{n-1} \frac{I_{\nu+k+1}(x)}{I_{\nu+k}(x)} +e^{-\frac{([x]-[\nu]-1)([x]+\nu+\nu_f)}{2x}}F([x]+\nu_f,x)(1+F([x]+\nu_f+1,x))\enskip.$$

Thus, it remains only to estimate the sum $ \sum_{n=1}^{[x]-[\nu]-1} \prod_{k=0}^{n-1} {I_{\nu+k+1}(x)}/{I_{\nu+k}(x)} $.  Using Corollary \ref{Hboundexp} again, we get

$$
\begin{aligned}
\prod_{k=0}^{n-1} \frac{I_{\nu+k+1}(x)}{I_{\nu+k}(x)}\leq \prod_{k=0}^{n-1} e^{-\alpha_0 \frac{\nu+k+\frac{1}{2}}{x}}&=\exp\left(-\frac{\alpha_0}{2x} n(2\nu+n)\right)=\exp\left(-\frac{\alpha_0}{2x} [(n+\nu)^2-\nu^2] \right) \enskip.
\end{aligned}
$$

Applying the same technique used for the lower bound, we have

\begin{equation}
\label{upper H}
\begin{aligned}
e^{\frac{1}{2x} \left(\nu+\frac{1}{2}\right)^2} \sum_{n=1}^{[x]-[\nu]-1} e^{-\frac{1}{2x} \left(n+\nu+\frac{1}{2}\right)^2}\leq  \sum_{n=1}^{[x]-[\nu]-1} \prod_{k=0}^{n-1} \frac{I_{\nu+k+1}(x)}{I_{\nu+k}(x)} \leq  e^{\frac{\alpha_0 \nu^2}{2x}} \sum_{n=1}^{[x]-[\nu]-1} e^{-\frac{\alpha_0}{2x} \left(n+\nu\right)^2}\enskip. \\
\end{aligned}
\end{equation}

Since both the upper and lower bounds are similar, we focus only on the upper bound.  The lower bound can be treated similarly.  
$$
\begin{aligned}
& \sum_{n=1}^{[x]-[\nu]-1} e^{-\frac{\alpha_0}{2x} \left(n+\nu\right)^2}\\
&=\sum_{k=[\nu]+1}^{[x]-1} e^{-\frac{\alpha_0}{2x} \left(k+\nu_f\right)^2}\\
&\leq \int_{[\nu]}^{[x]-1} e^{-\frac{\alpha_0}{2x} \left(y+\nu_f\right)^2} dy\\
&=\sqrt{\frac{2x}{\alpha_0}} \int_{\sqrt{\frac{\alpha_0}{2x}}\nu}^{\sqrt{\frac{\alpha_0}{2x}[x]-1+\nu_f}} e^{-u^2}du\\
&\leq \sqrt{\frac{2x}{\alpha_0}} \int_{\sqrt{\frac{\alpha_0}{2x}}\nu}^{\sqrt{\frac{\alpha_0}{2x}}x} e^{-u^2}du\\
&=\sqrt{\frac{2x}{\alpha_0}} \left[ \int_{\sqrt{\frac{\alpha_0}{2x}}\nu}^{\infty} e^{-u^2}du-\int_{\sqrt{\frac{\alpha_0}{2x}} x}^{\infty} e^{-u^2}du\right]\\
\end{aligned}
$$

\begin{equation}
\label{Gaussian Integral}
\leq \sqrt{\frac{2x}{\alpha_0}}\left[\frac{e^{-\frac{\alpha_0}{2x} \nu^2 }}{\sqrt{\frac{\alpha_0}{2x}} \nu +\sqrt{\frac{\alpha_0}{2x} \nu^2+\frac{4}{\pi}}}-\frac{e^{-\frac{\alpha_0}{2x} x^2}}{\sqrt{\frac{\alpha_0}{2x}} x+\sqrt{\frac{\alpha_0}{2x} x^2+2}}\right]
\end{equation}

where after a $u$-substitution, we have used the inequality (see \cite{AS}), 

$$\frac{e^{-x^2}}{x+\sqrt{x^2+2}} \leq \int_x^{\infty} e^{-t^2}dt \leq \frac{e^{-x^2}}{x+\sqrt{x^2+\frac{4}{\pi}}} \; \; \; x\geq 0.$$

From (\ref{upper H}) and (\ref{Gaussian Integral}), we have

\begin{equation}
\label{H5}
\begin{aligned}
&H(\nu,x) \leq \sqrt{\frac{2x}{\alpha_0}}\left[\frac{1}{\sqrt{\frac{\alpha_0}{2x}} \nu +\sqrt{\frac{\alpha_0}{2x} \nu^2+\frac{4}{\pi}}}-\frac{e^{-\frac{\alpha_0}{2x} (x^2-\nu^2)}}{\sqrt{\frac{\alpha_0}{2x}} x+\sqrt{\frac{\alpha_0}{2x} x^2+2}}\right] \\
&\hspace{1in}+e^{-\frac{\alpha_0}{2x} ([x]-[\nu]-1)([x]+\nu+\nu_f-1)}\frac{F\left([x]+\nu_f-\frac{1}{2}\right)}{1-F\left([x]+\nu_f+\frac{1}{2},x\right)}
\end{aligned}
\end{equation}

which can be written as

$$
\begin{aligned}
H(\nu,x)& \leq \frac{2x}{\alpha_0}\left[\frac{1}{\nu +\sqrt{\nu^2+\frac{8x}{\pi\alpha_0}}}-\frac{e^{-\frac{\alpha_0}{2x} (x^2-\nu^2)}}{x+\sqrt{x^2+\frac{4x}{\alpha_0}}}\right] \\
&\hspace{1in}+e^{-\frac{\alpha_0}{2x} ([x]-[\nu]-1)([x]+\nu+\nu_f-1)}\frac{F\left([x]+\nu_f-\frac{1}{2}\right)}{1-F\left([x]+\nu_f+\frac{1}{2},x\right)}\enskip ,
\end{aligned}
$$

yielding the upper bound in (\ref{Our Upper H Bound}).

To complete the proof then, we just need to prove the lower bound.  Repeating  similar arguments above,
$$
\begin{aligned}
\sum_{n=1}^{[x]-[\nu]-1} e^{-\frac{1}{2x} \left(n+\nu+\frac{1}{2}\right)^2} &\geq \int_{[\nu]+1}^{[x]-\nu+[\nu]+1} e^{-\frac{1}{2x}\left(y+\frac{1}{2}\right)^2}dy\\
&=\sqrt{2x}\int_{\frac{[\nu]+\frac{3}{2}}{\sqrt{2x}}}^{\frac{[x]-\nu+[\nu]+\frac{3}{2}}{\sqrt{2x}}} e^{-u^2}du\\
&\geq \sqrt{2x}\int_{\frac{\nu+\frac{3}{2}}{\sqrt{2x}}}^{\frac{x-\nu_f+\frac{3}{2}}{\sqrt{2x}}}e^{-u^2}du\\
&=\sqrt{2x} \left[ \int_{\frac{1}{\sqrt{2x}}\left(\nu+\frac{3}{2}\right)}^{\infty} e^{-u^2}du -  \int_{\frac{1}{\sqrt{2x}}\left([x]-\nu_f+\frac{3}{2}\right)}^{\infty} e^{-u^2}du\right]\\
&\geq \sqrt{2x}\left[ \frac{e^{-\frac{1}{2x} \left(\nu+\frac{3}{2}\right)^2}}{\frac{1}{\sqrt{2x}} \left(\nu+\frac{3}{2}\right)+\sqrt{\frac{1}{2x} \left(\nu+\frac{3}{2}\right)^2+2}}\right.\\
&\left.\hspace{1in}- \frac{e^{-\frac{1}{2x} \left([x]-\nu_f+\frac{3}{2}\right)^2}}{\frac{1}{\sqrt{2x}} \left([x]-\nu_f+\frac{3}{2}\right)+\sqrt{\frac{1}{2x} \left([x]-\nu_f+\frac{3}{2}\right)^2+\frac{4}{\pi}}}\right] \enskip.
\end{aligned}
$$
\end{enumerate}

Thus,

$$
\begin{aligned}
&H(\nu,x)\geq e^{\frac{1}{2x} \left(\nu+\frac{1}{2}\right)^2}\sqrt{2x}\left[ \frac{e^{-\frac{1}{2x} \left(\nu+\frac{3}{2}\right)^2}}{\frac{1}{\sqrt{2x}} \left(\nu+\frac{3}{2}\right)+\sqrt{\frac{1}{2x} \left(\nu+\frac{3}{2}\right)^2+2}}- \frac{e^{-\frac{1}{2x} \left([x]-\nu_f+\frac{3}{2}\right)^2}}{\frac{1}{\sqrt{2x}} \left([x]-\nu_f+\frac{3}{2}\right)+\sqrt{\frac{1}{2x} \left([x]-\nu_f+\frac{3}{2}\right)^2+\frac{4}{\pi}}}\right]\\
&\hspace{3in}+e^{-\frac{([x]-[\nu]-1)([x]+\nu+\nu_f)}{2x}}F([x]+\nu_f,x)(1+F([x]+\nu_f+1,x)) \\
&\Rightarrow H(\nu,x)\geq e^{-\frac{1}{x}\left(\nu+1\right)} \frac{\sqrt{2x}}{\frac{1}{\sqrt{2x}} \left(\nu+\frac{3}{2}\right)+\sqrt{\frac{1}{2x} \left(\nu+\frac{3}{2}\right)^2+2}}\\
&\hspace{1in}-e^{-\frac{1}{2x} ([x]-\nu-\nu_f+1)([x]+\nu-\nu_f+2)} \frac{\sqrt{2x}}{\frac{1}{\sqrt{2x}} \left([x]+\frac{3}{2}\right)+\sqrt{\frac{1}{2x} \left([x]+\frac{3}{2}\right)^2+\frac{4}{\pi}}}\\
&\hspace{1in}+e^{-\frac{([x]-[\nu]-1)([x]+\nu+\nu_f)}{2x}}F([x]+\nu_f,x)(1+F([x]+\nu_f+1,x))\\
\end{aligned}
$$

which is the same as

$$
\begin{aligned}
H(\nu,x)&\geq   \frac{2xe^{-\frac{1}{x}\left(\nu+1\right)}}{\nu+\frac{3}{2}+\sqrt{\left(\nu+\frac{3}{2}\right)^2+4x}}- \frac{2xe^{-\frac{1}{2x} ([x]-\nu-\nu_f+1)([x]+\nu-\nu_f+2)}}{[x]+\frac{3}{2}+\sqrt{ \left([x]+\frac{3}{2}\right)^2+\frac{8x}{\pi}}}\\
&\hspace{2in}+e^{-\frac{([x]-[\nu]-1)([x]+\nu+\nu_f)}{2x}}F([x]+\nu_f,x)(1+F([x]+\nu_f+1,x))\enskip.
\end{aligned}
$$

Theorem \ref{Hboundexp} is proved.

\begin{flushright}
$\square$
\end{flushright}

\begin{figure}[htp!]
\includegraphics[width=17cm, height=10cm]{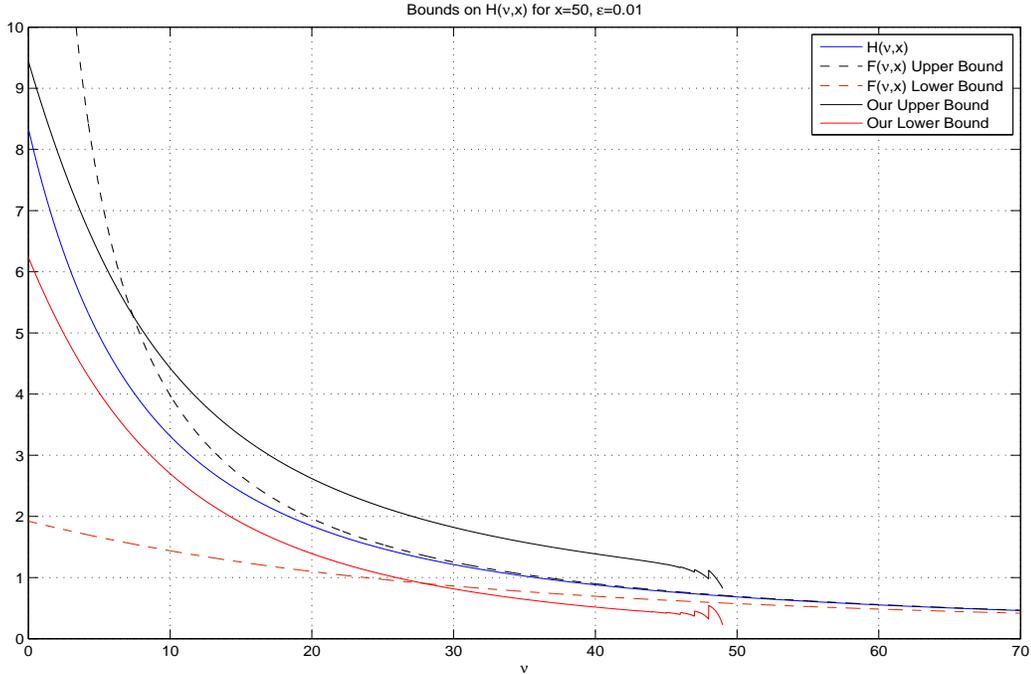}
\caption{A comparison of our bounds (\ref{Our H Bounds}), (\ref{Our Lower H Bound}) and (\ref{Our Upper H Bound}) on $H(\nu,x)$ compared to the true value of $H(\nu,x)$ for $x=50$ for $\nu \in [0,200]$ taken in steps of $0.01$.  The $F(\nu,x)$ Lower/Upper bounds refer to those in (\ref{Geometric Bounds}).  The value $\epsilon=0.01$ is chosen to truncate the infinite sum of Bessel functions occurring in the numerator of $H(\nu,x)$ so that the terms beyond a certain index are less than $\epsilon$.  We note that near $\nu=0$, our bounds are substantially better than using (\ref{Geometric Bounds}) and scale like $\sqrt{x}$ for large $x$ which we would not have been able to obtain otherwise, and is the main purpose of this paper.}  
\label{Hboundplot}
\end{figure}

We illustrate the bounds (\ref{Our H Bounds}), (\ref{Our Lower H Bound}) and (\ref{Our Upper H Bound}) on $H(\nu,x)$ in Figure  \ref{Hboundplot} for $x=50$ and $\nu\in [0,70]$ in steps of $0.01$ and also plot the true values of $H(\nu,x)$ computed using MATLAB.  For comparison, we also plot the $F(\nu,x)$ Lower/Upper bounds (\ref{Geometric Bounds}).  The value $\epsilon=0.01$ is chosen to truncate the infinite sum of Bessel functions occurring in the numerator of $H(\nu,x)$ so that the terms beyond a certain index are less than $\epsilon$.   We notice that there are regimes in $\nu$ for which our lower and upper bounds are worse and better than those using the geometric series-type bounds, but that near $\nu=0$, our bounds are substantially better, and is a result of the first difference in (\ref{Our Lower H Bound}) and (\ref{Our Upper H Bound}) obtained by the use of the exponential approximations (\ref{Bessel Ratio I}) on $I_{\nu+1}(x)/I_{\nu}(x)$.

\section{Main Results: Applications}
\label{appls}

In this section, we give some applications of Theorem \ref{Hboundexp}.  First, we briefly review the Skellam distribution, and relate it to the function $H(\nu,x)$.

  Let $X_1\sim Pois(\lambda_1)$ and $X_2\sim Pois(\lambda_2)$ be two independent Poisson random variables with parameters $\lambda_1$ and $\lambda_2$, respectively.  Then, the distribution of the random variable $W=X_1-X_2$ is called a Skellam distribution with parameters $\lambda_1$ and $\lambda_2$.  We denote this by $W=X_1-X_2\sim Skellam(\lambda_1,\lambda_2)$ and have,

$$\mathbb{P}\left[W=n\right]=e^{-\left(\lambda_1+\lambda_2\right)} \left( \frac{\lambda_1}{\lambda_2}\right)^{\frac{n}{2}} I_{|n|}\left(2\sqrt{\lambda_1\lambda_2}\right).$$

The probabilistic value of $H(\nu,x)$ is now immediate: if $\lambda_1=\lambda_2=\lambda>0$, then $H(\nu,2\lambda)={\mathbb{P}\left[W>\nu\right]}/{\mathbb{P}\left[W=\nu\right]}$.  The quantity

$$\mathcal{H}(\nu,2\lambda)=\frac{1}{H(\nu,2\lambda)+1}=\frac{\mathbb{P}\left[W=\nu\right]}{\mathbb{P}\left[W\geq \nu\right]}$$ is important in the actuarial sciences for describing the probability of death at time $\nu$ given death occurs no earlier than time $\nu$, and is known as the hazard function.

\subsection{Application 1:  Bounds on $\exp(-x)I_{\nu}(x)$ for $x\geq 0$ and $\nu\in \mathbb{N}$ and the Skellam$(\lambda_1,\lambda_2)$ Mass Function}
\label{Asymptotic application}

Since the distribution of $W\sim Skellam(\lambda,\lambda)$ is symmetric, we have

\begin{equation}
\label{Skellam Symmetry}
\exp(-2\lambda)I_0(2\lambda)=\mathbb{P}\left[W=0\right]=\frac{1}{2H(0,2\lambda)+1}.
\end{equation}

Thus, we may apply the bounds on  ${I_{\nu+1}(x)}/{I_{\nu}(x)}$ and $H(\nu,x)$  given in Corollary \ref{Bessel Ratio} and Theorem \ref{Hboundexp}, respectively, to obtain sharp upper and lower bounds on $\exp(-x)I_0(x)$ and hence on $$\exp(-x)I_{\nu}(x)=\prod_{k=0}^{\nu-1} \frac{I_{k+1}(x)}{I_k(x)} \exp(-x)I_0(x)$$ for $\nu\in \mathbb{N}$.  We note that this result therefore improves the asymptotic formula $$\exp(-x)I_{\nu}(x)\sim \frac{1}{\sqrt{2\pi x}}  \; \; \; {\rm as} \; \; \; x\rightarrow \infty$$ but only for $\nu\in \mathbb{N}$, and in particular, gives a bound on $\mathbb{P}\left[W=\nu\right]$ for $W\sim Skellam(\lambda,\lambda)$ by setting $x=2\lambda$.

\begin{Theorem}
\label{x Large}
Set  $\alpha_0=-\log(\sqrt{2}-1)\approx 0.8814.$  Then, for $\nu\in \mathbb{N}$ and $x\geq 0$,

\begin{enumerate}

\item If $\nu\leq x$, $$\frac{\exp((-\frac{\nu^2}{2x}\left(\frac{\nu+1}{\nu}\right))}{1+2\mathcal{U}(0,x)} \leq \exp(-x)I_{\nu}(x)\leq \frac{\exp(-\frac{\alpha_0}{2x} \nu^2)}{1+2\mathcal{L}(0,x)}$$

\item  If $\nu> x$, $$\frac{e^{-\frac{[x]^2}{2x}\left( \frac{[x]+1}{[x]}\right)}B\left([x]+\frac{x}{2}+1,\nu-[x]\right)  \frac{(x/2)^{\nu-[x]}}{(\nu-[x]-1)!}}{1+2\mathcal{U}(0,x)}\leq e^{-x}I_{\nu}(x)\leq \frac{e^{-\frac{\alpha_0}{2x}[x]^2} B\left([x]+x+\frac{1}{2},\nu-[x]\right) \frac{x^{\nu-[x]}}{(\nu-[x]-1)!}}{{1+2\mathcal{L}(0,x)}}
$$

where $B(x,y)$ denotes the Beta function and $\mathcal{L}(\nu,x)$ and $\mathcal{U}(\nu,x)$ are the lower and upper bounds, respectively, from Theorem \ref{Hboundexp}.
\end{enumerate}

\end{Theorem}

{\bf Proof of Theorem \ref{x Large}:}  

\begin{enumerate}

\item  By Corollary \ref{Bessel Ratio}, for $k+1\leq x$,

$$
e^{-\frac{k+1}{x}}\leq \frac{I_{k+1}(x)}{I_{k}(x)} \leq e^{-\alpha_0 \frac{k+\frac{1}{2}}{x}}.
$$

so that for $\nu\leq x$,

$$e^{-\frac{\nu^2}{2x} \left(\frac{\nu+1}{\nu}\right) }= \prod_{k=0}^{\nu-1} e^{-\frac{k+1}{x}} \leq \prod_{k=0}^{\nu-1} \frac{I_{k+1}(x)}{I_k(x)}\leq \prod_{k=0}^{\nu-1} e^{-\alpha_0 \frac{k+\frac{1}{2}}{x}}=e^{-\frac{\alpha_0}{2x} \nu^2}\enskip.$$

Thus,

$$e^{-\frac{\nu^2}{2x}\left(\frac{\nu+1}{\nu}\right)}e^{-x}I_0(x)\leq e^{-x}I_{\nu}(x)\leq e^{-\alpha_0 \frac{k+\frac{1}{2}}{x}}=e^{-\frac{\alpha_0}{2x} \nu^2}e^{-x}I_0(x)$$

since

$$e^{-x}I_{\nu}(x)=\prod_{k=0}^{\nu-1} \frac{I_{k+1}(x)}{I_k(x)} e^{-x}I_0(x) .$$

By (\ref{Skellam Symmetry}) then, we have for $\nu\leq x$,

$$\frac{e^{-\frac{\nu^2}{2x}\left(\frac{\nu+1}{\nu}\right)}}{1+2\mathcal{U}(0,x)} \leq e^{-x}I_{\nu}(x)\leq \frac{e^{-\frac{\alpha_0}{2x} \nu^2}}{1+2\mathcal{L}(0,x)} \enskip.$$

\item To prove the second assertion in theorem \ref{x Large}, notice that for $\nu> x$, $$\prod_{k=0}^{\nu-1} \frac{I_{k+1}(x)}{I_k(x)} = \prod_{k=0}^{[x]-1} \frac{I_{k+1}(x)}{I_k(x)} \prod_{k=[x]}^{\nu-1} \frac{I_{k+1}(x)}{I_k(x)}$$ and each term in the first product has $k\leq x$ so that by the previous argument,

$$e^{-\frac{[x]^2}{2x}\left( \frac{[x]+1}{[x]}\right)}\prod_{k=[x]}^{\nu-1} \frac{I_{k+1}(x)}{I_k(x)}\leq \prod_{k=0}^{\nu-1} \frac{I_{k+1}(x)}{I_k(x)} \leq e^{-\frac{\alpha_0}{2x}[x]^2} \prod_{k=[x]}^{\nu-1} \frac{I_{k+1}(x)}{I_k(x)} .$$

Next,

$$
\begin{aligned}
\frac{I_{k+1}(x)}{I_k(x)} &\leq \frac{x}{k+\frac{1}{2} + \sqrt{\left(k+\frac{1}{2}\right)^2 + x^2}} \leq \frac{x}{k+\frac{1}{2}+x}
\end{aligned}
$$

so that

$$
\begin{aligned}
\prod_{k=[x]}^{\nu-1} \frac{I_{k+1}(x)}{I_k(x)} \leq \prod_{k=[x]}^{\nu-1} \frac{x}{k+\frac{1}{2}+x}&=\frac{x^{\nu-[x]} \Gamma\left([x]+x+\frac{1}{2}\right)}{\Gamma\left(\nu+\frac{1}{2}+x\right)}\\
&=B\left([x]+x+\frac{1}{2},\nu-[x]\right) \frac{x^{\nu-[x]}}{(\nu-[x]-1)!}
\end{aligned}
$$

and similarly, using $\sqrt{a^2+b^2}\leq a+b$ for $a,b\geq 0$,

$$
\begin{aligned}
\prod_{k=[x]}^{\nu-1} \frac{I_{k+1}(x)}{I_k(x)} &\geq \prod_{k=[x]}^{\nu-1} \frac{x}{k+1 + \sqrt{\left(k+1\right)^2 + x^2}} \\
&\geq  \prod_{k=[x]}^{\nu-1} \frac{x}{2(k+1)+x}\\
&=  \prod_{k=[x]}^{\nu-1} \frac{x/2}{k+1+x/2}\\
&=\frac{(x/2)^{\nu-[x]} \Gamma\left([x]+1+\frac{x}{2}\right)}{\Gamma\left(\nu+\frac{x}{2}+1\right)}\\
&=B\left([x]+\frac{x}{2}+1,\nu-[x]\right)  \frac{(x/2)^{\nu-[x]}}{(\nu-[x]-1)!}\\
\end{aligned}
$$

Thus, since $\exp(-x)I_{\nu}(x)=\prod_{k=0}^{\nu-1} {I_{k+1}(x)}/{I_k(x)} \exp(-x)I_0(x)$, we get

$$\frac{e^{-\frac{[x]^2}{2x}\left( \frac{[x]+1}{[x]}\right)}B\left([x]+\frac{x}{2}+1,\nu-[x]\right)  \frac{(x/2)^{\nu-[x]}}{(\nu-[x]-1)!}}{1+2\mathcal{U}(0,x)}\leq e^{-x}I_{\nu}(x)\leq \frac{e^{-\frac{\alpha_0}{2x}[x]^2} B\left([x]+x+\frac{1}{2},\nu-[x]\right) \frac{x^{\nu-[x]}}{(\nu-[x]-1)!}}{{1+2\mathcal{L}(0,x)}}
.$$

Thus theorem \ref{x Large} is proved.

\end{enumerate}
\begin{flushright}
$\square$
\end{flushright}

A few remarks of Theorem \ref{x Large} are in order:

\begin{enumerate}
\item   One may simplify the upper bound using the bounds found in \cite{Alzer} on the Beta function, $B(x,y)$,

$$\alpha\leq \frac{1}{xy} -B(x,y)\leq \beta$$

where $\alpha=0$ and $\beta=0.08731\ldots$ are the best possible bounds.

\item By setting $x=2\sqrt{\lambda_1\lambda_2}$ and multiplying (1) and (2) in Theorem \ref{x Large} through by $$\left(\sqrt{\frac{\lambda_1}{\lambda_2}}\right)^{\nu}\exp\left[-\left(\sqrt{\lambda_1}+\sqrt{\lambda_2}\right)^2\right],$$ we obtain precise bounds on $\mathbb{P}\left[W=\nu\right]$ for $W\sim Skellam(\lambda_1,\lambda_2)$.

\item  It's important to note that if one were to use the geometric series-type bound (\ref{Geometric Bounds}) with $\nu=0$, that one would not achieve the behavior of $1/\sqrt{x}$ that we have in Theorem \ref{x Large} which is indeed, guaranteed by the asymptotic $\exp(-x)I_{\nu}(x)\rightarrow 1/\sqrt{2\pi x}$ as $x\rightarrow \infty$ and exhibited by our non-asymptotic bounds on $H(\nu,x)$ in (\ref{Our Lower H Bound}) and (\ref{Our Upper H Bound}).

\end{enumerate}

\begin{figure}[htp!]
\includegraphics[width=15cm, height=10cm]{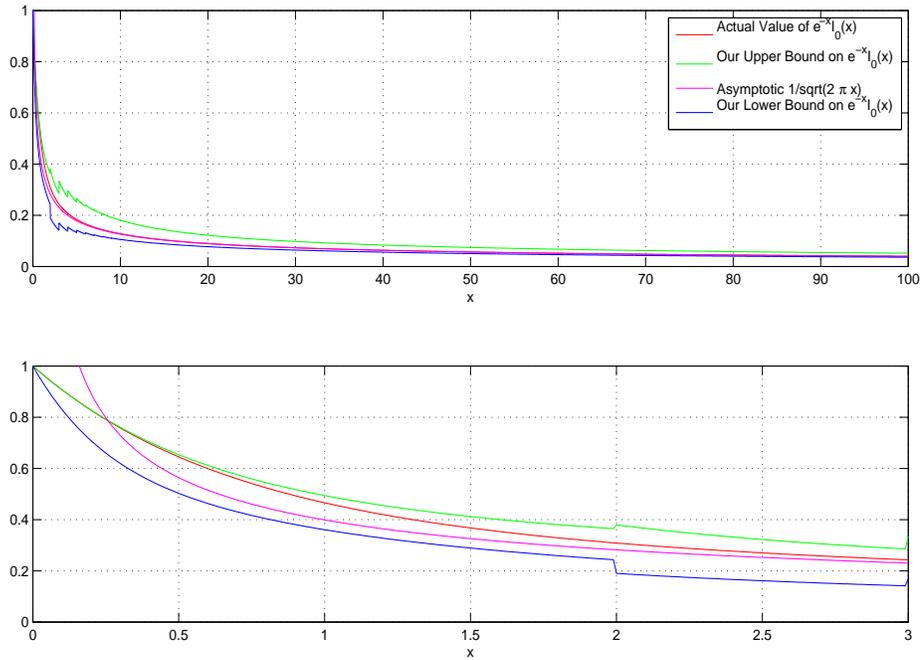}
\caption{A comparison of our bounds on $\exp(-x)I_0(x)$ to the true value and the asymptotic $1/\sqrt{2\pi x}$.  All plots are taken in steps of $1/100$ on a window of $[0,100]$.  Top panel:  behavior for $x\in [0,100]$.  Bottom panel:  the same plot, but over the interval $x\in [0,3]$.  We note that for $[x]<2$, our upper bound is much more accurate than the asymptotic $1/\sqrt{2\pi x}$ and that in general is quite good.  For large $x$, all curves $\sim 1/\sqrt{x}$ in agreement with the asymptotic $1/\sqrt{2\pi x}$, a behavior that would not have been achieved if instead one used the geometric series-type bounds (\ref{Geometric Bounds}) in place of $\mathcal{L}(0,x)$ and $\mathcal{U}(0,x)$.}
\label{besseli0}
\end{figure}

As an example of applying our bounds non-asymptotically, we plot $\exp(-x)I_0(x)$, its asymptotic ${1}/{\sqrt{2\pi x}}$ as $x\rightarrow \infty$, and the functions ${1}/({2\mathcal{L}(0,x)+1})$ and ${1}/(2\mathcal{U}(0,x)+1)$ in Figure \ref{besseli0}.  In steps of $1/100$, over the interval $[0,100]$, the top panel illustrates the behavior of all these functions over the interval $[0,100]$.  We note that for large values of $x$, all functions values converge to zero, are extremely close, and are all on the order of $1/\sqrt{x}$ - something that one would not see by using instead the naive geometric-type bounds (\ref{Geometric Bounds}) with $\nu=0$.   In the second panel, we restrict to the interval $[0,3]$ as our bounds transition across the line $[x]=2$.  We note that for $[x]<2$, our upper bound is much more accurate than the asymptotic $1/\sqrt{2\pi x}$ and that in general is quite good.

\subsection{Application 2:  Concentration Inequality for the Skellam Distribution}
\label{Skellam application}
We now present a concentration inequality for the $Skellam(\lambda,\lambda)$, the proof of which is a direct consequence of Theorem \ref{x Large} and the identity,

$$\mathbb{P}\left[|W|>\nu\right] = 1-\mathbb{P}\left[-\nu\leq W \leq  \nu\right] =2H(\nu,2\lambda) \exp(-2\lambda)I_{\nu}(2\lambda) .$$

\begin{cor}
Let $W\sim Skellam(\lambda, \lambda)$, and define $\alpha_0=-\log\left(\sqrt{2}-1\right)$.  Then,

\begin{enumerate}

\item If $\nu\leq x$, $$\frac{2\exp\left(-\frac{\nu^2}{4\lambda}\left(\frac{\nu+1}{\nu}\right)\right)}{1+2\mathcal{U}(0,2\lambda)} \leq \frac{\mathbb{P}\left[|W|>\nu\right]}{H(\nu,2\lambda)} \leq \frac{2\exp\left(-\frac{\alpha_0}{4\lambda} \nu^2\right)}{1+2\mathcal{L}(0,2\lambda)}$$

\item  If $\nu>x$, $$
\begin{aligned}
\frac{2\exp\left(-\frac{[2\lambda]^2}{4\lambda}\left( \frac{[2\lambda]+1}{[2\lambda]}\right)\right)B\left([2\lambda]+\lambda +1,\nu-[2\lambda]\right)  \frac{\lambda ^{\nu-[2\lambda]}}{(\nu-[2\lambda]-1)!}}{1+2\mathcal{U}(0,2\lambda)}&\leq \frac{\mathbb{P}\left[|W|>\nu\right]}{H(\nu,2\lambda)}\\
\frac{2\exp\left(-\frac{\alpha_0}{4\lambda}[2\lambda]^2\right) B\left([2\lambda]+2\lambda+\frac{1}{2},\nu-[2\lambda]\right) \frac{(2\lambda)^{\nu-[2\lambda]}}{(\nu-[2\lambda]-1)!}}{{1+2\mathcal{L}(0,2\lambda)}}&\geq  \frac{\mathbb{P}\left[|W|>\nu\right]}{H(\nu,2\lambda)}\\
\end{aligned}
$$

where $B(x,y)$ denotes the Beta function and $\mathcal{L}(\nu,x)$ and $\mathcal{U}(\nu,x)$ are the lower and upper bounds, respectively, from Theorem \ref{Hboundexp}.
\end{enumerate}
\end{cor}

\section{Summary and Conclusions}
\label{conclusion}

In \cite{Viles}, the function $H(\nu,x)=\sum_{n=1}^{\infty} {I_{\nu+n}(x)}/{I_{\nu}(x)}$ for $x\geq 0$ and $\nu\in \mathbb{N}$ appears as a key quantity in approximating a sum of dependent random variables that appear in statistical estimation of network motifs as a Skellam$(\lambda,\lambda)$ distribution .  A necessary scaling of $H(\nu,x)$ at $\nu=0$ of $\sqrt{x}$ is necessary, however, in order for the error bound of the approximating distribution to remain finite for large $x$.  In this paper, we have presented a quantitative analysis of $H(\nu,x)$ for $x,\nu\geq 0$ necessary for these needs in the form of upper and lower bounds in Theorem \ref{Hboundexp}.   Our technique relies on bounding current estimates on ${I_{\nu+1}(x)}/{I_{\nu}(x)}$ from above and below by quantities with nicer algebraic properties, namely exponentials, while optimizing the rates when $\nu+1\leq x$ to maintain their precision.

In conjunction with the mass normalizing property of the $Skellam(\lambda,\lambda)$ distribution, we also give applications of this function in determining explicit error bounds, valid for any $x\geq 0$ and $\nu\in \mathbb{N}$, on the asymptotic approximation $\exp({-x})I_{\nu}(x)\sim {1}/{\sqrt{2\pi x}}$ as $x\rightarrow \infty$, and use them to provide precise upper and lower bounds on $\mathbb{P}\left[W=\nu\right]$ for $W\sim Skellam(\lambda_1,\lambda_2)$.  In a similar manner, we derive a concentration inequality for the $Skellam(\lambda,\lambda)$ distribution, bounding $\mathbb{P}\left[|W|\geq \nu\right]$ where $W\sim Skellam(\lambda,\lambda)$ from above and below.

While we analyze the function $H(\nu,x)$ $\nu\in \mathbb{N}$, $x\geq 0$ for our purposes, we leave as future research the analysis for non integer $\nu$, as well as consideration of the generalized function

$$H(\nu,x)=\sum_{n=1}^{\infty}\left(\frac{\lambda_1}{\lambda_2}\right)^{\frac{n}{2}} \frac{  I_{\nu+n}(2\sqrt{\lambda_1\lambda_2})}{I_{\nu}(2\sqrt{\lambda_1\lambda_2})}$$
that would appear for the $Skellam(\lambda_1,\lambda_2)$ distribution.  We hope that the results laid here will form the foundation of such future research in this area.  

It is also unknown as to whether normalization conditions for $\{\exp({-x})I_{\nu}(x)\}_{\nu=-\infty}^{\infty}$ induced by the $Skellam(\lambda,\lambda)$  hold for $\nu$ in a generalized lattice $\mathbb{N}+\alpha$, and if so, what the normalizing constant is.  Such information would provide a key in providing error bounds on the asymptotic $\exp({-x})I_{\nu}(x)\sim {1}/{\sqrt{2\pi x}}$ for non-integer values of $\nu$.

\bibliographystyle{elsarticle-num}

\begin{thebibliography}{00}
\footnotesize

\bibitem{AS}
{\sc M. Abramowitz and I.A. Stegun,} Handbook of Mathematical Functions with Formulas, Graphs, and Mathematical Tables, Dover, New York, 1972.  

\bibitem{Amos}
{\sc D. E. Amos,}  Computation of Modified Bessel Functions and Their Ratios, Math. Comp. 28 (1974) 239-251.

\bibitem{Alzer}
{\sc H. Alzer.}  Sharp Inequalities for the Beta Function, Indagationes Mathematicae 12, 1 (2001) 15-21.

\bibitem{Barbour}
{\sc A.D. Barbour and L. Chen}, An Introduction to Stein's Method, Singapore University Press, Singapore, 2005.

\bibitem{Baricz Turanians}
{\sc A. Baricz,} Bounds for Turanians of Modified Bessel Functions, arXiv:1202.4853 (2013).

\bibitem{Baricz Edin}
{\sc A. Baricz,} Bounds for Modified Bessel Functions of the First and Second Kinds, Proceedings of the Edinburgh Mathematical Society 52 (2010) 575-599.

\bibitem{Baricz Neuman}
{\sc A. Baricz and E. Neuman,} Inequalities Involving Modified Bessel functions of the First Kind,   Journal of Mathematical Analysis and Applications 332 (2007) 265-271.

\bibitem{Baricz Sum}
{\sc A. Baricz and T. Pogany}, On a Sum of Modifed Bessel Functions,  arXiv:1301.5429 (2013).

\bibitem{Baricz Marcum Q}
{\sc A. Baricz and Y. Sun,} New bounds for the Generalized Marcum Q-function, IEEE Transactions on Information Theory 55 (2009), 7.

\bibitem{Cranston}
{\sc M.C. Cranston and S.A. Molchanov,} On a concentration inequality for sums of independent isotropic vectors, Electron. Commun. Probab. 17(27) (2012), 1-8.

\bibitem{Fotopoulos}
{\sc S. B. Fotopoulos and K.J. Venkata,} Bessel inequalities with applications to conditional log returns under GIG scale mixtures of normal vectors, Stat. Probab. Lett. 66 (2004) 117-125.

\bibitem{Hwang}
{\sc Y. Hwang,} Difference-Based Image Noise Modeling Using Skellam Distribution, Pattern Analysis and Machine Intelligence 34, 7 (2012) 1329-1341.

\bibitem{Johnson}
{\sc N.L. Johnson,} On an Extension of the Connection between Poisson and $\chi^2$-distributions, Biometrika 46 (1959) 352-363.

\bibitem{Kanter}
{\sc M. Kanter,} Probability inequalities for convex sets and multidimensional concentration, Journal of Multivariate Analysis 6, 2 (1976) 222-236.

\bibitem{Karlis}
{\sc D. Karlis and I. Ntzoufras,} Analysis of sports data using bivariate Poisson models, Journal of the Royal Statistical Society Series D 52 (3) (2003) 381-393.

\bibitem{Laforgia}
{\sc A. Laforgia and P. Natalini,} Some Inequalities for Modified Bessel Functions, Journal of Inequalities and Applications (2010) 1-10.

\bibitem{Luke}
{\sc Y. Luke,} Inequalities for Generalized Hypergeometric Functions, Journal of Approximation Theory 5 (1972) 41-65.

\bibitem{Marchand1}
{\sc E. Marchand and F. Perron,} Improving on the MLE of a bounded normal mean, Ann. Statist.  29 (2001) 1078-1093.

\bibitem{Marchand2}
{\sc E. Marchand and F. Perron,} On the minimax estimator of a bounded normal mean, Stat. Probab. Lett. 58 (2002) 327-333.

\bibitem{Neuman}
{\sc E. Neuman,} Inequalities Involving Modified Bessel Functions of the First Kind, Journal of Mathematical Analysis and Applications 171 (1992) 532-536.

\bibitem{Robert}
{\sc C. Robert,} Modified Bessel functions and their applications in probability and statistics, Stat. Probab. Lett. 9 (1990) 155-161.

\bibitem{Skellam}
{\sc J. Skellam,} The frequency distribution of the difference between two Poisson variates belonging to different populations, Journal of the Royal Statistical Society Series A 109 (3) (1946) 296.

\bibitem{Simpson}
{\sc H. Simpson and S. Spector,} Some monotonicity results for ratios of modified Bessel functions, Journal of Inequalities and Applications 42, 1 (1984) 95-98.

\bibitem{Simon}
{\sc M. Simon and M.S. Alouini,}  Digital Communication Over Fadding Channels: A Unified Approach to Performance Analysis, Wiley, New York, 2000.

\bibitem{Viles}
{\sc W. Viles, P. Balachandran and E. Kolaczyk.}  A Central Limit Theorem for Network Motifs.  Manuscript.

\bibitem{Wolfe}
{\sc P. Wolfe and K. Hirakawa,}  Efficient Multivariate Skellam Shrinkage for Denoising Photon-Limited Image Data: An Empirical Bayes Approach, Proc. IEEE Int. Conf. Image Processing (ICIP-09), Cairo, Egypt, Nov. 7-11 (2009)  pp. 2961-2964.

\bibitem{Yuan}
{\sc L. Yuan and J.D. Kalbfleisch,} On the Bessel distribution and related problems, Ann. Inst. Statist. Math. 52(3) (2000) 438-447.



\end{thebibliography}

\end{document}